\documentstyle[12pt]{amsart}
\newtheorem{prop}{Proposition}
\newtheorem{th}[prop]{Theorem}
\newtheorem{cor}[prop]{Corollary}
\newtheorem{lem}[prop]{Lemma}

\theoremstyle{definition}
\newtheorem{ack}{Acknowledgments} 
\theoremstyle{remark}
\newtheorem{rem}{Remark}

\newcommand{\bbR}{{\Bbb{R}}}
\newcommand{\calP}{{\cal{P}}}

\newcommand{\al}{\alpha}
\newcommand{\be}{\beta}

\newcommand{\de}{\delta}
\newcommand{\e}{\varepsilon}
\newcommand{\tht}{\theta}

\newcommand{\k}{\kappa}
\newcommand{\la}{\lambda}
\newcommand{\m}{\mu}
\newcommand{\n}{\nu}

\newcommand{\s}{\sigma}

\newcommand{\h}{\chi}

\newcommand{\Om}{\Omega}

\renewcommand{\span}{\operatorname{span}}
\newcommand{\supp}{\operatorname{supp}}
\newcommand{\card}{\operatorname{card}}
\newcommand{\codim}{\operatorname{codim}}

\newcommand{\sgn}{\operatorname{sgn}}

\newcommand{\disp}{\displaystyle}
\newcommand{\lb}{\label}

\newcommand{\emp}{\emptyset}

\newcommand{\lra}{\longrightarrow}
\newcommand{\str}{\stackrel}
\newcommand{\wtw}{if and only if }
\newcommand{\buo}{without loss of generality }

\makeatletter
\def\@currentlabel{2.1}\label{e:dispaa}
\def\@currentlabel{2.21}\label{e:dispau}
\def\@currentlabel{2.22}\label{e:dispav}
 \def\@currentlabel{2.23}\label{e:dispaw}
\def\@currentlabel{2.24}\label{e:dispax}
\makeatother

\makeatletter
\def\alphenumi{%
  \def\theenumi{\alph{enumi}}%
  \def\p@enumi{\theenumi}%
  \def\labelenumi{(\@alph\c@enumi)}}
\makeatother

\begin{document}

\title{Contractive projections and isometries in sequence spaces}
\author{Beata Randrianantoanina}
\address{Department of Mathematics \\ Bowling Green State University \\
         Bowling Green, OH 43403}
\email{brandri@@andy.bgsu.edu}
\maketitle

\begin{abstract}
We characterize 1-complemented subspaces of finite codimension in
strictly monotone
one-$p$-convex, $2<p<\infty,$ sequence spaces. Next we describe, up
to isometric isomorphism, all possible types of 1-unconditional
structures in  sequence spaces with few surjective isometries. We also
give a new example of a class of real sequence spaces with few
surjective isometries.
\end{abstract}

\section{Introduction} \lb{flinnin}

This paper is divided into three parts. Throughout we consider real
sequence spaces with 1-unconditional basis.

First we study images of contractive projections --
we prove (Theorem~\ref{pro}) that in strictly monotone and
one-$p$-convex, $2<p<\infty,$ (or, dually, one-$q$-concave, $1<q<2$)
sequence spaces
every 1-complemented subspace of finite codimension $n$
contains all but at most $2n$
basic vectors.
Calvert and Fitzpatrick \cite{CF86} showed that if any such
 hyperplane  is
1-complemented then the space is isometric to  $\ell_p$ or $c_0.$

  Characterizations of contractive projections are important in
  approximation
   theory and there exists an extensive literature on the subject
   (see
\cite{ChP} and \cite{BP88}   for the detailed
disscusion and references).

Theorem~\ref{pro} applies to a  rich class of spaces including e.g.
 $\ell_p,\ 1<p<\infty,\ p \neq 2,$
$\ell_p(\ell_r)$ where $2<p,r<\infty,$ or  $1<p,r<2,$ as well as
 a wide
class of Orlicz and Lorentz spaces. It generalizes the analogous
result known for classical sequence spaces: see \cite{BCh,Wph89,Wph92} for
 $\ell_1$, \cite{BM,BP88} for
$\ell_p,\ 1<p<\infty ,\ p \neq 2,$ \cite{Wph} for
$\ell_p^n ,\ 1 \le p<\infty ,\ p\neq 2.$
The analogous result is not true in $c_0$ \cite{BCh} or
$\ell^\infty$ \cite{B88}.

Our method of proof is quite different, and we believe simpler, than
those used before.

Next we investigate all (up to isometric equivalence) 1-unconditional
structures in a given sequence space. This is an isometric version of
the question of uniqueness of unconditional basis, which has been
studied
since late sixties (c.f.~\cite{BCLT} for various sequence spaces and
\cite{K93} for detailed disscussion and references).

In the complex case the situation is well understood.
 Kalton and Wood proved \cite[Theorem 6.1]{KW} that all 1-unconditional
 bases in a complex Banach space are isometrically equivalent (cf. also
 \cite[discussion on page 452 and Corollary 3.13]{R86} and \cite{FJ74}).
  Lacey and Wojtaszczyk
\cite{LW} observed that this does not hold in real $L_p$-spaces -- they
give a complete description of the two possible types of 1-unconditional
  structure
in $L_p$ (cf. also \cite{BL77}).
As far as we know very little work has been done since then in real
Banach spaces (except \cite{R86}).

In Theorem~\ref{gsiso} below we establish that in real sequence spaces
which have few surjective isometries there are two  types of isometrically
non-equivalent 1-unconditional
  structure.
 Corollary~\ref{2siso} formulates
some additional assumptions which yield the
uniqueness of 1-unconditional basis.

It now becomes of interest to describe the spaces satisfying assumptions of
Theorem~\ref{gsiso} -- i.e. spaces with few surjective isometries. This is
a problem that have been studied for its own right by many authors starting
with Banach \cite{Ban} who characterized isometries in $\ell_p^n.$ In the
complex case the theory is  well developed (see e.g. the survey
\cite{FJs} and its references).

In the real case Braveman and Semenov \cite{BS} (cf. also
\cite[Theorem 9.8.3]{Rol}) proved that symmetric sequence spaces have few
(in our sense) isometries. Skorik \cite{Sk} showed an analogous result
for
a special
class of
real sequence spaces.
We do not know of any other pertinent references.

In Section~\ref{iseq} we provide a new example of a class of spaces with
only elementary surjective isometries.
As an application of
Theorem~\ref{pro}, we prove
(Theorem~\ref{siso})  that all
surjective isometries between two strictly monotone
sequence spaces which are both one-$p$-convex, $2<p<\infty,$ or
one-$q$-concave,
$1<q<2,$
 are elementary. Our results
are valid in both finite and infinite-dimensional spaces.

\begin{ack}
I wish to express my gratitude to Professor Nigel Kalton for his interest
in this work and many valuable disscussions.
\end{ack}

\section{Norm-one complemented subspaces of finite codimension in
sequence
spaces}
\lb{pseq}

We say that a Banach space $X$ is {\bf one-$\bold p$-convex} (resp.
{\bf one-$\bold q$-concave})  \lb{defpcon}
 if
for every
choice of elements $\{ x_i \}_{i=1}^n $ in $X$ the following inequality
holds:
\begin{align*}
   \| \left( \sum_{i=1}^n |x_i|^p \right)^{1/p} \| \ \ &\le \ \ \left(
\sum_{i=1}^n
\|x_i\|^p
\right)^{1/p} \ \ \ \ \ \ \ \hfill \text{if } 1\le p < \infty, \\
\intertext{or, respectively,}
   \| \left( \sum_{i=1}^n |x_i|^q \right)^{1/q} \| \ \ &\ge \ \ \left(
\sum_{i=1}^n
\|x_i\|^q
\right)^{1/q} \ \ \ \ \ \ \ \hfill \text{if } 1\le q < \infty,
\end{align*}
(cf. \cite[Definition 1.d.3]{LT2}).

\begin{th} \lb{pro}
Let  $X$ be a
strictly monotone
 sequence space $(\dim X = d \ge 3)$ with a 1-unconditional basis
$\{e_i \}_{i=1}^d.$ Suppose that

$(a)$ \  $X$  is one-$p$-convex, $2<p<\infty,$

\noindent or

$(b)$ \  $X$ is
one-$q$-concave, $1<q<2,$  and smooth at each basic vector.

Then any
1-complemented subspace $F$ of codimension $n$ in $X$ contains all
 but at most $2n$ basic
vectors of $X.$
\end{th}

\begin{rem}
Notice that Theorem~\ref{pro} states only necessary and not sufficient
conditions for the subspace to be 1-complemented (unlike the theorem of
Baronti and Papini \cite{BP88} for $\ell_p$). Also Baronti and Papini
\cite{BP88} prove that in $\ell_p$ every 1-complemented subspace of
finite codimension is an intersection of 1-complemented hyperplanes.
The analogous statement is not true in general (cf.
\cite{BG}).
\end{rem}

For the proof of Theorem~\ref{pro} we will need the
 following observation which we state in the form
of the lemma for easy reference
\begin{lem} \lb{lpro}
Let $X$ be a one-$p$-convex, $2<p<\infty,$ sequence space with a
1-unconditional basis
 and let
$P:X \str{\text{onto}}{\lra} F$ be a projection. Assume that there exist
disjoint elements $x,y \in X$ such that $\supp Py \supset
\supp x, \ Px = x$ and $\card(\supp x) < \infty.$

Then $\|P\| > 1.$
\end{lem}

\begin{pf*}{Proof of Lemma~\ref{lpro}}
Let us assume, for contradiction, that $\|P\| \le 1$ and take $x, y$ with
$\|x\| = \|y\| = 1.$
By one-$p$-convexity of $X$ and  since $x$ and $y$ are disjoint we
get for all $\tau \in
\bbR$:
\begin{equation} \lb{pro1}
 \| P(x + \tau y)\| \le \| x + \tau y \| = \|(|x|^p + |\tau y|^p)^{1/p}
 \|
\le (1 + |\tau |^p)^{1/p}
\end{equation}

Since $p>2$ $X$ is one-2-convex (\cite[Proposition 1.d.5, p.49]{LT2})  and
for any
$\tau \in
\bbR$ we get:
\begin{equation} \lb{pro2}
\begin{split}
\| \left(|P(x + \tau y)|^2 + |P(x - \tau y)|^2\right)^{1/2} \| \
   &\le \left(\|P(x + \tau e_1)\|^2 + \|P(x - \tau e_1)\|^2\right)^{1/2} 
\\
   &\underset{\text{by (\ref{pro1})}}{\le} \sqrt{2}\,(1+
|\tau|^p)^{1/p} \end{split}
\end{equation}

On the other hand
\begin{equation} \lb{pro3}
\begin{split}
\| (  |P(x &+ \tau y)|^2 + |P(x - \tau y)|^2
)^{1/2}
\|
\\
  &\ge \| \sum_{i\in\supp x} \left( {2}\,\left(|x_i|^2 +
\tau^2 |(Py)_i|^2\right)^2\right)^{\frac12}\, e_i \| \\
  &\ge\sqrt{2}\,\| \sum_{i\in\supp x} |x_i| \sqrt{1 +
\tau^2\left(\frac{(Py)_i}{x_i}\right)^2} \,e_i \| \\
  &\ge \sqrt{2}\, \sqrt{1 + \eta \tau^2}\/\|x\| =
\sqrt{2}\, \sqrt{1 + \eta \tau^2}
\end{split}
\end{equation}
where ${\disp \eta = \min_{i\in\supp x}
\left\{\left(\frac{(Py)_i}{x_i}\right)^2\right\}.}$ Notice that $\eta >0,$
since $\supp x \subset\supp Py.$

Combining (\ref{pro2}) and (\ref{pro3}) we get
$\sqrt{1 + \eta \tau^2} \le  (1+ |\tau|^p)^{1/p} $ which gives us the
desired contradiction when $|\tau | < \eta^{\frac{1}{p-2}}.$
\end{pf*}

\begin{pf*}{Proof of Theorem~\ref{pro}}
We first prove part $(a)$ of the theorem.
Let  $F$ be a 1-complemented subspace of codimension n, say
${\disp F = \bigcap_{j=1}^n \ker f_j}$ for some $f_j \in X^*$
and the contractive projection $P\lra F$ is given by
 ${\disp P = Id_X - \sum_{j=1}^n f_j
\otimes v^j}$ for some linearly independent  $v^j \in X$ with
$f_j(v^k) = \de_{jk}$ (where
$\de_{jk}$ denotes Kronecker delta).

Assume that $e_i \not\in F$ if $i \in I.$ If $\card I < n$ there is
nothing to prove so \buo
$\{1,2,\dots,n\} \subset I$ and $f_j(e_i) = \de_{ij}, \ i,j \le n.$

Notice first that if ${\disp l\not\in \bigcup_{i=1}^n \supp v^i}$
then ${\disp P(e_l) = e_l  - \sum_{i=1}^n f_i(e_l) v^i}$ and so
$\left(P(e_l)\right)_l = 1.$ Thus, by strict monotonicity of $X$,
$P(e_l) = e_l,$ i.e. $e_l \in F.$ Therefore
\begin{equation} \lb{pro-1}
I \subset \bigcup_{i=1}^n \supp v^i.
\end{equation}

Now take any  ${\disp a = \sum_{i=1}^n a_i e_i.}$
Then
\begin{equation*}
\begin{split}
      P(a) &= a - \sum_{j=1}^n f_j(a) v^j
               = \sum_{i=1}^n a_i e_i  -
                           \sum_{j=1}^n \sum_{i=1}^n a_i f_j(e_i) v^j\\
                   &= \sum_{i=1}^n a_i e_i - \sum_{i=1}^n a_i v^i.
\end{split}
\end{equation*}

Hence there exists $a_0 \in \span\{e_1,\dots,e_n\}$ such that
$${\disp \supp P(a_0) \setminus \{1,\dots,n\} =
\bigcup_{i=1}^n \supp v^i \setminus \{1,\dots,n\}.}$$

If ${\disp \card\left(\bigcup_{i=1}^n \supp v^i \setminus
 \{1,\dots,n\}\right) > n}$ then
$\card(\supp P(a_0) \setminus \{1,2,\dots,n\} ) \ge n+1$ and there
exists
$ x\in F$ with $\supp x \subset \supp P(a_0)
\setminus \{1,2,\dots,n\} $ (since $\codim F = n < n +1$).
Now $x$ and $a_0$ satisfy assumptions of
 Lemma~\ref{lpro} which contradicts the fact that $P$ is contractive.

 Thus ${\disp \card\left(\bigcup_{i=1}^n \supp v^i \setminus
 \{1,\dots,n\}\right) \le n}$ and, by~\eqref{pro-1},
 $$\card I \le \card\left(\bigcup_{i=1}^n \supp v^i \right) \le 2n$$
 which proves part $(a)$ of the theorem.

We prove part $(b)$ by duality.
 Consider contractive projection
 ${\disp P^* = Id_{X^*} - \sum_{j=1}^n v^j\otimes f_j }.$
$X^*$ is one-$p$-convex for some $p>2$ and strictly monotone so
  by
part $(a)$ we get that, say,
$v^1,\dots,v^n \subset \span\{e_1^*,\dots,e_{2n}^*\}.$ Thus,
by~\eqref{pro-1},
$I \subset \{1,\dots,2n\}$ since $X$ is strictly monotone.
\end{pf*}

  Theorem~\ref{pro} can be  combined with our previous results
   about nonexistence
  of 1-com\-ple\-men\-ted hyperplanes in nonatomic function spaces
   which do not have
  any bands isometrically equal to $L_2$ \cite[Theorem 2]{Ran1}
   (cf. also
\cite[Theorem 4.3]{KR}, \cite[Theorem 2.7]{th}).

\begin{cor}
Suppose that $X$ is a separable
strictly monotone  function space on $(\Om,\m)$ which is either
one-$p$-con\-vex for some $2<p<\infty,$ or one-$q$-concave for some
$1<q<2$ and smooth at $\h_A$ for every atom $A$ of $\m.$ Suppose further
that for some $g \in X^*,\ $ $\ker g$ is 1-complemented in $X.$ Then $g$
is of the form $\al \h_A + \be \h_B,$ where $\al, \be \in \bbR$ and
$A,B$ are atoms of $\m.$
\end{cor}

The above  statement  exactly parallels (and extends) the
theorem proved
by Beauzamy and Maurey for $L_p$ \cite[Proposition 3.1, p.135]{BM}.

\section{Isometries and 1-unconditional bases of sequence spaces}
\lb{bases}

An operator $T:X \lra Y$ between two sequence spaces with
1-un\-con\-di\-tion\-al
bases $\{e_i\}_{i=1}^d$ and $\{f_i\}_{i=1}^d,$ resp. ($d\le\infty$), will
be called {\bf elementary} if
$$T(e_i) = a_i f_{\s(i)}$$
for some $a_i \in \bbR$ and a permutation $\s$ of $\{1,\dots,d\}.$

We will say that a pair of indices $k,l$ is {\bf interchangeable} in $X$
if for any $x,z \in X$ \ $|x_k| = |z_l|, |x_l| = |z_k| $ and $|x_i| =
|z_i|$
for all $i \neq k,l$ imply that $\|x\| = \|z\|.$ Space $X$ is rearrangement
invariant \wtw every two indices are interchangeable.

\begin{th} \lb{gsiso}
Suppose that $X,Y$  are separable sequence spaces with
1-un\-con\-di\-tion\-al
bases $\{e_i\}_{i=1}^d$ and $\{f_i\}_{i=1}^d,$ resp. and suppose that all
surjective isometries of one of the spaces  $X$ or $Y$ onto itself are
elementary.
Suppose that $T:X \lra Y$ is a surjective isometry.

Then there exist a set $A \subset \{1,\dots,d\}$ and a 1-1 map $\s :
A \lra \s(A)
\subset \{1,\dots,d\}$ such that for every $i \in A$
$$T(e_i) = \e_i f_{\s(i)}$$
where $\e_i = \pm 1.$

The complementary sets $B_X = \{1,\dots,d\} \setminus A$ and $B_Y =
\{1,\dots,d\}
\setminus \s(A)$ split into families of disjoint pairs $\cal P_X
\subset 2^{B_X},
\cal P_Y \subset 2^{B_Y}$ so that there exists a 1-1 map $\tau :
\cal P_X \str{
\text{onto}}{\lra} \cal P_Y$ and if $\tau(i,j) = (k,l)$ then
\begin{align*}
T(e_i) &= \frac{\de_i}{\|f_k + f_l \|} (f_k + \e_i f_l) \\
T(e_j) &= \frac{\de_i}{\|f_k + f_l \|} (f_k - \e_i f_l)
\end{align*}
where $\de_i, \e_i = \pm 1.$

Moreover \newline $(a)$ all pairs $(i,j) \in \cal P_X$ and $(k,l) \in
\cal P_Y$ are interchangeable in $X$ or $Y,$ resp.,
\newline $(b)$ if all isometries of $Y$ (resp. $X$) onto itself are
elementary then the set $A$ (resp. $\s(A)$) depends only on the spaces
$X, Y$ and not on the isometry $T.$
\end{th}

The following fact is an immediate consequence of Theorem~\ref{gsiso}.

\begin{cor} \lb{2siso}
In the situation of  Theorem~\ref{gsiso} if we assume additionally that
no 2-dimensional subspace of one of the spaces $X$ or $Y$ is isometric to
$\ell_2^2$ and both spaces $X$ and $Y$ are either one-2-convex or
one-2-concave then every surjective isometry $T:X \lra Y$ is elementary.
\end{cor}

\begin{rem}
Since all surjective isometries of rearrangement-invariant sequence spaces
onto itself are elementary \cite{BS} Corollary~\ref{2siso} may be viewed
as an isometric and sequence space version of the deep result of Kalton
about (isomorphic) uniqueness of lattice structure in nonatomic
2-convex
(or strictly 2-concave) Banach lattices which embed complementably in a
strictly 2-convex
(resp. strictly 2-concave) r.i. function space
\cite[Theorems~8.1 and 8.2]{K93}.
\end{rem}

{\it Proof of Theorem~\ref{gsiso}.}
We use
all the notation as introduced above.

Let us first see that the final remark follows readily from the
main
statment of the theorem.

$(a)$ Let $y \in Y$ and $x = T^{-1}(y).$ Consider the element $\tilde{x}
\in X$
such that $\tilde{x}_j = - x_j$ and $\tilde{x}_{\n} = x_{\n}$ for $\n
\neq j.$ Then
$\| \tilde{x} \| = \|x\|$ and so $\|y\| = \|T\tilde{x}\|.$ But from the
form of $T$
we see that $(T\tilde{x})_k = \e_iy_l, \ (T\tilde{x})_l = \e_iy_k$ and
$(T\tilde{x})_{\n}
= y_{\n}$ for $\n \neq k,l.$ Hence $(k,l)$ is interchangeable in $Y.$
Proof for $(i,j) \in \calP_X$ is similar.

$(b)$  Assume that the set $A$ depends on the isometry $T$ and let use
the notation $A_T$ to emphisize that dependence. Assume that $i \in A_U
\setminus A_T$ for some isometries $U, T.$ Then
$$ TU^{-1}(\e_if_{\s(i)}) = T(e_i) =
\frac{\de_i}{\|f_k + f_l \|} (f_k + \e_i f_l)  $$
which contradicts the fact that the isometry $ TU^{-1} : Y \lra Y$ is
elementary.

Now let us return to the proof of the main statement of the theorem.
It is  clear that if $T$ has the described form
then so does $T^{-1}.$ So we can assume \buo that all isometries of, say,
$Y$
onto itself are elementary.

We will split the proof of the theorem into a series of lemmas.

\begin{lem} \lb{lgsiso1}
\mbox{      }\nopagebreak

\noindent $(a)$ For any $i \le d$ there exist at most two indices $k,l$
such that
$$T(e_i) = \al_k f_k + \al_l f_l.$$

\noindent   $(b)$ If $k \neq l$ then $\al_k,\al_l \neq 0.$
\end{lem}

\begin{lem} \lb{lgsiso2}
Suppose that for some $i,j,k,l \le d $
\begin{align*}
T(e_i) &= \al_k f_k + \al_l f_l \\
T(e_j) &= \be_k f_k + \be_m f_m,
\end{align*}
where $\al_k,\be_k \neq 0.$ Then

\noindent $(a)$ $l = m, \ \al_l,\be_m \neq 0$ and $\sgn(\al_k
\al_l)
=-
\sgn(
\be_k \be_m).$ \newline
   $(b)$ $|\al_k| = |\al_l| = |\be_k| = |\be_m|.$
\end{lem}

\begin{lem} \lb{lgsiso3}
Suppose that for all $n \le d$  $\card\supp T(e_n) \le 2.$ Let $i,k,l \le d
$
be such that
$$T(e_i) = \al_k f_k + \al_l f_l,$$
where $\al_k,\al_l \neq 0, k \neq l.$ Then there exist a unique $j \neq i,$
and $\be_k, \be_l \neq 0$ so that
$$T(e_j) = \be_k f_k + \be_l f_l.$$
\end{lem}

\begin{pf*}{Proof of Lemma~\ref{lgsiso1}}
Denote
\begin{align*}
T(e_j) &= \sum_{m=1}^d \al_{j,m}f_m \\
T^{-1}(f_n) &= \sum_{j=1}^d \be_{n,j}e_j
\end{align*}
For any choice of signs $\e = (\e_j)_{j=1}^d, \ \e_j = \pm 1, $
we define an isometry $S_{\e}:X \lra X$ by $S_{\e}(e_j) =
\e_je_j.$ By unconditional convergence we get for every $n$:
\begin{equation*}
\begin{split}
 TS_{\e}T^{-1}(f_n)&= T\left(\sum_{j=1}^d \be_{n,j}\e_je_j\right) =
                  \sum_{j=1}^d \e_j \be_{n,j}\left(
                  \sum_{m=1}^d \al_{j,m}f_m \right)   \\
                &= \sum_{m=1}^d \left(\sum_{j=1}^d \e_j\be_{n,j}
                  \al_{j,m}\right)f_m.
\end{split}
\end{equation*}

Since $TS_{\e}T^{-1}$ is elementary (as a surjective isometry of $Y$)
we conclude that for
every $ n \le d$ and $\e = (\e_j)_{j=1}^d$ there exists exactly
one $m$ such that
\begin{equation}  \lb{gsiso1}
\sum_{j=1}^d \e_j\be_{n,j}
                  \al_{j,m}    \neq 0.
\end{equation}

Now fix $ i \le d.$ Since $T^{-1}$ is onto there exists $ n \le d$ with
$\be_{n,i} \neq 0.$ By (\ref{gsiso1}) we get:
\begin{align}
\exists !\ k \text{ with  } \ \ \ \
&\sum_{j=1}^d \be_{n,j}
                  \al_{j,k}    \neq 0
\ \ \ \ \ \ \ \ \ \ \ \ \  \ \ \ \ \ \ \ \
(\e_j = 1 \text{ for all }j) \lb{gsiso2} \\
\exists !\ l \text{ with  }  \ \
-&\sum_{j \neq i} \be_{n,j}
                  \al_{j,l} + \be_{n,i}\al_{i,l}   \neq 0\  \ \ \ \ \ \ \
\left( \e_j = \begin{cases}
                -1& \ j \neq i \\
                 1& \ j=i
              \end{cases}
\ \ \right) \lb{gsiso3}
\end{align}

Hence  for any
 $m \neq k,l$ $\al_{i,m} = 0$ i.e.
$$T(e_i) = \al_{i,k} f_k + \al_{i,l} f_l$$
which proves part $(a)$ of the lemma.

  Part $(b)$ follows immediately from
equations (\ref{gsiso2}) and (\ref{gsiso3}).
\end{pf*}

\begin{pf*}{Proof of Lemma~\ref{lgsiso2}}
\mbox{   }

\noindent  $(a)$ \  Let $c= -\al_k \be_k^{-1}.$ Then since $T$ is an
isometry we have
$$\|\al_lf_l + c\be_mf_m\| = \|T(e_i + ce_j)\| = \|T(e_i - ce_j)\| =
\|2\al_kf_k + \al_lf_l - c\be_mf_m\|.$$

Hence $l=m, \ \al_l,\be_m \neq 0$ and $\sgn\al_l = -\sgn(c\be_m)$ i.e.
$\sgn(\al_k \al_l) = -\sgn(\be_k \be_m).$

\noindent  $(b)$ \ Since $l =m$ denote $T(e_j) = \be_k f_k + \be_l f_l.$
Then
\begin{align*}
T^{-1}(f_k) &= -\be_lMe_i + \al_lMe_j \\
T^{-1}(f_l) &= \be_kMe_i - \al_kMe_j,
\end{align*}
where $M = (\al_l\be_k - \al_k \be_l)^{-1}$ ($\al_l\be_k - \al_k \be_l
\neq 0$ by part $(a)$).

Denote by $S$ the isometry of $X$ such that $S(e_i) = -e_i$ and
$S(e_j) = e_j.$ Then
$$TST^{-1}(f_k) = T(\be_lMe_i + \al_lMe_j)
                  = M(\al_l\be_k + \al_k \be_l)f_k + 2M\be_l\al_lf_l.$$
Since $TST^{-1}$ is a surjective isometry of $Y$ it is elementary
and since $2M\be_l\al_l \neq 0$ we get $\al_l\be_k + \al_k \be_l =0$
i.e.
\begin{equation} \lb{lgsiso21}
\al_l\be_k = - \al_k \be_l.
\end{equation}
Moreover
\begin{equation} \lb{lgsiso22}
|2M\be_l\al_l| = 1.
\end{equation}
Combining \eqref{lgsiso21} and \eqref{lgsiso22} we obtain
$|\al_k| = |\al_l|$ and $ |\be_k| = |\be_l|,$ and since
$\|T(e_i)\| = \|T(e_j)\|$ we have $|\al_k| = |\al_l| = |\be_k| = |\be_l|.$
\end{pf*}

\begin{pf*}{Proof of Lemma~\ref{lgsiso3}}
Lemma~\ref{lgsiso2} implies that for any $j \neq i$ we have either
$\supp T(e_j) = \supp T(e_i)$ or $\supp T(e_j) \cap \supp T(e_i) = \emp.$
Hence, by surjectivity of T, there exists $j \neq i$ with
$T(e_j) = \be_k f_k + \be_l f_l$ and by Lemma~\ref{lgsiso2} $\be_k,
\be_l \neq 0.$

Uniqueness of $j$ is an immediate consequence of the fact that $T$
preserves the dimension of subspaces.
\end{pf*}

\section{Isometries in one-$p$-convex sequence spaces} \lb{iseq}

\begin{th} \lb{siso}
Suppose $X,\ Y$ are separable strictly monotone sequence spaces with
1-un\-con\-di\-tional bases
 and $\dim X = \dim Y = d \ge 3 \ (d \le \infty ).$ Suppose that

$(a)$ \  $X,Y$  are one-$p$-convex, $2<p<\infty,$

\noindent or

$(b)$ \  $X,Y$ are
one-$q$-concave, $1<q<2,$  and smooth at each basic vector.

 Then
any isometry $U$ from $X$ onto $Y$ is of the form
$$U \left( \sum_{k=1}^d a_ke_k \right) = \sum_{k=1}^d
\e_ka_kf_{\sigma(k)} $$
where $\s$ is a permutation of $\{1,\dots
,d \}$ and $\e_k = \pm 1$ for $k=1,\dots
,d.$
\end{th}

\begin{pf}
We will prove the theorem with the assumption $(a)$. Part $(b)$ follows
by duality.

For any $k \le d$ the hyperplane
$\{ x_k = 0 \}$
 is 1-complemented in $X$
and so is $U\{ x_k = 0 \}$ in $Y.$ By Theorem~\ref{pro} there are at
most two numbers $k_1,k_2 \le d $ such that $U\{ x_k = 0 \} = \{
\al_1y_{k_1} + \al_2y_{k_2} = 0 \}$ for some $\al_1, \al_2 \in \bbR.$ We
will say that coordinates $k,l$ are related if $\{k_1,k_2\} \cap
\{l_1,l_2\} \neq \emp.$

For the proof of the theorem we need three technical lemmas.

\begin{lem} \lb{lsiso1}
Suppose
$U\{ x_k = 0 \} = \{
\al_1y_{k_1} + \al_2y_{k_2} = 0 \}$ where  $\al_1, \al_2 \neq 0$ and
suppose that $l$ is related to $k.$ Then $\{l_1,l_2\} \subset
\{k_1,k_2\}.$
\end{lem}

\begin{lem} \lb{lsiso2}
For any $k \le d$ there is at most one coordinate $l$ $( \neq k)$
related to $k.$
\end{lem}

\begin{lem} \lb{lsiso3}
For any $k \le d$ there exist $i,j \le d, \ \k_i,\k_j \in \bbR$ such that
$$U(e_k) = \k_if_i + \k_jf_j.$$

Moreover if both
$\k_i,\k_j \neq 0$ then there exist (unique) $l
\neq k$ and  $\la_i,\la_j \in \bbR$ such that
$U(e_l) = \la_if_i + \la_jf_j.$
\end{lem}

Let us first see that Theorem~\ref{siso} indeed follows from
Lemma~\ref{lsiso3}.

If, say, $\k_j = 0$ then $|\k_i| = 1$ since $U$ is an isometry. So we
need only to show that $\k_i,\k_j$ cannot both be nonzero.

Assume, for contradiction, that  $\k_i,\k_j \neq 0.$ Then by
Lemma~\ref{lsiso3} there exists $l \neq k$ such that
$U(e_l) = \la_if_i + \la_jf_j$ for some $\la_i,\la_j \in \bbR.$
By one-$p$-convexity of $Y$
 we get:
\begin{align*}
1 &= \| \k_if_i + \k_jf_j \|  \le
 (\k_i^p + \k_j^p)^{1/p}
 < (\k_i^2 + \k_j^2)^{1/2} \\
1 &= \| \la_if_i + \la_jf_j \|
  \le  (\la_i^p + \la_j^p)^{1/p}
  \le  (\la_i^2 + \la_j^2)^{1/2}
\end{align*}

Hence\ \ \  \ \ \ \ \ \ \ \ \ \ \ \ \ \  \ \ \ \ \ \ \ \ \ \  $ \k_i^2 +
\k_j^2
+
\la_i^2
+
\la_j^2
>
2.$

So, say,
\begin{equation} \lb{siso1}
 (\k_i^2 + \la_i^2)^{1/2} = \|(\k_i,\la_i)\|_2 > 1
\end{equation}
On the other hand by one-2-convexity of $X$ for any $(a,b) \in \bbR^2$ we
have $\|ae_k + be_l \| \le \|(a,b)\|_2.$ But
\begin{equation*}
\begin{split}
  \|ae_k + be_l \| &= \| (a\k_i + b\la_i)f_i + (a\k_j + b\la_j)f_j \| \\
                   &\ge \| (a\k_i + b\la_i)f_i\| = | (a\k_i + b\la_i)|
\end{split}
\end{equation*}
So $$ | (a\k_i + b\la_i)| \le \|(a,b)\|_2 $$ and this means that
$\|(\k_i,\la_i)\|_2 \le 1$ which contradicts (\ref{siso1}) and ends the
proof of the theorem.
\end{pf}

\begin{pf*}{Proof of Lemma~\ref{lsiso1}}
Our assumption is
\begin{equation} \lb{1}
U\{ x_k = 0 \} = \{
\al_1y_{k_1} + \al_2y_{k_2} = 0 \}
\end{equation}
where  $\al_1, \al_2 \neq 0, $ and $l$ is related to $k.$
Without loss of generality  $l_1 = k_1$ and we have
\begin{equation} \lb{2}
U\{ x_l = 0 \} = \{
\be_1y_{k_1} + \be_2y_{k_2} = 0 \}
\end{equation}
where  $\be_1 \neq 0. $ If $\be_2 = 0$ there is nothing to prove so let
us assume $\be_2 \neq 0.$
Proposition~\ref{pro} applied to the isometry $U^{-1}$ gives us:
\begin{align}
U\{ y_{k_1} = 0 \} &= \{
\m_1x_{m_1} + \m_2x_{m_2} = 0 \} = H_{k_1}  \lb{3} \\
U\{ y_{k_2} = 0 \} &= \{
\n_1x_{n_1} + \n_2x_{n_2} = 0 \} = H_{k_n}  \lb{4} \\
U\{ y_{l_2} = 0 \} &= \{
\tht_1x_{t_1} + \tht_2x_{t_2} = 0 \} = H_{l_1}  \lb{5}
\end{align}
Denote \ \ $E_l = \{x_l = 0 \}, E_k = \{x_k = 0 \} \subset X.$

Since $U^{-1}$ is an isometry equations (\ref{1}), (\ref{3}), (\ref{4})
imply that $E_k \cap H_{k_1} =E_k \cap H_{k_2} =H_{k_1}  \cap H_{k_2}$
i.e. the following systems of equations are equivalent:
$$
\left\{
\begin{matrix}
        \m_1x_{m_1} + \m_2x_{m_2} = 0 \\ x_k = 0
\end{matrix} \right\}
\equiv
\left\{
\begin{matrix}
        \n_1x_{n_1}\ +\ \n_2x_{n_2} = 0 \\ x_k = 0
\end{matrix} \right\}
\equiv
\left\{
\begin{matrix}
        \m_1x_{m_1} + \m_2x_{m_2} = 0 \\ \n_1x_{n_1}\ +\ \n_2x_{n_2} = 0
\end{matrix} \right\}
$$

Since these systems have rank 2 this implies that, say, $m_1 = n_1 =k,\
m_2 = n_2 \neq k$ and $\m_2,\n_2 \neq 0.$

Similarly by considering equations (\ref{2}), (\ref{3}), (\ref{5})
 we
obtain $m_1 = t_1,\ m_2 =t_2$ and either $m_1 =l$ or $m_2 = l.$ Hence
 $k = m_1 = n_1 = t_1$ and $l = m_2 = n_2 =t_2.$
But this
means that $\codim( H_{k_1} \cap H_{k_2}  \cap H_{l_2}) \le 2.$ Since
$U$
is an isometry we have $$\codim\{y_{k_1},y_{k_2},y_{l_2} = 0 \} = \codim
U( H_{k_1} \cap H_{k_2}  \cap H_{l_2})  = \codim
( H_{k_1} \cap H_{k_2}  \cap H_{l_2})  \le 2.$$
Hence $l_2 = k_2.$
\end{pf*}

\begin{pf*}{Proof of Lemma~\ref{lsiso2}}
If $k$ is related to $l$ then for at least one of $k,l,$ say $k,$
$U\{x_k=0\} = \{\al_1y_{k_1} + \al_2y_{k_2} = 0 \} $ where $\al_1, \al_2
\neq 0.$ Then by Lemma~\ref{lsiso1} $\{l_1,l_2\} \subset \{k_1,k_2\},$ so
if $t$ is related to $l$ it is also related to $k$ and
$\{t_1,t_2\} \subset \{k_1,k_2\}.$ But then $U\{x_k,x_l,x_t = 0 \}
\subset \{ y_{k_1},y_{k_2} = 0 \}$ and so $t \in \{k,l\}.$
\end{pf*}

\begin{pf*}{Proof of Lemma~\ref{lsiso3}}
We have $${\disp U(e_k) \in \bigcap_{\n \neq k} U\{x_{\n} = 0 \}}.$$

By Lemma~\ref{lsiso2} there exists at most one coordinate $l$ related to
$k$ and by Lemma~\ref{lsiso1} $\{k_1,k_2,l_1,l_2\} = \{i,j\} $ where $i
\neq j$ \wtw there exists $l \neq k$ related to $k.$ Moreover
\begin{equation} \lb{lsiso3.1}
\text{if }\ \  \n \neq k,l\ \ \ \ \ \  \text{ then }\ \ \ \ \  \{\n_1,\n_2\}
\cap
\{i,j\}
= \emptyset
\end{equation}
where $U\{x_\n = 0 \} = \{\al(\n)y_{\n_1} + \be(\n)y_{\n_2} = 0 \}.$

Since $U$ is 1-1 and onto
$$\bigcap_{\n \neq k,l} U\{x_{\n} = 0 \} = \bigcap_{\m \neq i,j} \{
y_{\m} = 0 \}.$$
Hence $U(e_k),U(e_l) \in \span\{f_i,f_j\}$ which proves the first part of
the lemma.

The second part follows immediately from the fact that $\span \{
U(e_k),U(e_l) \} = \span\{f_i,f_j\}$ and condition~\eqref{lsiso3.1}.
\end{pf*}

 \bibliographystyle{standard}
    \bibliography{tref}

\end{document}